\documentclass[11pt]{amsart}
\usepackage{upref,amsmath,amssymb,graphicx}

\newtheorem{theorem}{Theorem}[section]

\newtheorem*{lemma*}{Lemma}

\theoremstyle{definition}

\theoremstyle{remark}

\numberwithin{equation}{section}

\newcommand{\abs}[1]{\lvert#1\rvert}

\newcommand{\dist}{\operatorname{dist}}
\begin{document}

\title[Conformal contractions]{Conformal contractions and lower bounds on the density of harmonic measure}
\author{Leonid V. Kovalev}
\address{215 Carnegie, Mathematics Department, Syracuse University, Syracuse, NY 13244}
\email{lvkovale@syr.edu}
\thanks{Supported by the National Science Foundation grant DMS-1362453.}

\maketitle

\begin{abstract}
We give a concrete sufficient condition for a simply-connected domain to be the image of the unit disk under a nonexpansive conformal map. This class of domains is also characterized by having sufficiently dense harmonic measure. The relation with the harmonic measure provides a natural higher-dimensional analogue of this problem, which is also addressed.  
\end{abstract}

\section{Introduction}\label{intro}

The images of the unit disk $\mathbb D$ under conformal maps $f$ with the normalization $f(0)=0=f'(0)-1$ have long been understood and characterized in terms of their Green's function, capacity of the complement, and so on (e.g., the books \cite{Go} and \cite{Ts} expose this circle of ideas). This paper studies the effect of a uniform bound on the derivative of a conformal map: namely, $|f'(z)|\le 1$ for all $z\in\mathbb D$. This condition can be equivalently stated as  $|f(z)-f(w)|\le |z-w|$ for all $z,w\in\mathbb D$; such $f$ may be called a \textit{conformal contraction}. Under the normalization $f(0)=0$, it follows that the image $f(\mathbb D)$ must be contained in $\mathbb D$. However, not every subdomain of the unit disk is its image under a conformal contraction. 

Let us consider a convex domain $\Omega\subset \mathbb C$ that contains $0$ and has $C^{1,1}$-smooth boundary. With such a domain we associate three radii: 
\begin{itemize}
\item \emph{outer radius} $R_O$  is the smallest radius of a disk centered at $0$ and containing $\Omega$; 
\item \emph{inner radius} $R_I$  is the largest radius of a disk centered at $0$ and contained in $\Omega$; 
\item \emph{curvature radius} $R_C$ is the minimal radius of curvature of $\partial \Omega$. It is the largest radius $R$ such that $\Omega$ can be written as a union of open disks of radius $R$. 
\end{itemize} 

Note that $R_O\ge R_I$ and $R_O\ge R_C$, while there is no general relation between $R_I$ 
and $R_C$. 

\begin{theorem}\label{contra} Let $\Omega\subset \mathbb C$ be a convex domain that contains $0$ and has $C^{1,1}$-smooth boundary. If the radii $R_O$, $R_I$, and $R_C$ satisfy
\begin{equation}\label{sufc}
(R_O-R_C)\frac{\log R_I-\log R_C}{R_I-R_C}  + \frac12 \log R_C\le 0
\end{equation}
then $\Omega=f(\mathbb D)$ for some conformal map $f$ such that $f(0)=0$ and 
$\sup \abs{f'}\le 1$. (When $R_I=R_C$, the difference quotient is understood as $1/R_I$.)
\end{theorem} 

We will also consider the harmonic measure of domain $\Omega$ with respect to $0$, denoted $\omega_\Omega(\cdot, 0)$. In the context of Theorem~\ref{contra}, of particular interest is the Radon-Nikodym derivative of $\omega_\Omega(\cdot, 0)$ with respect to arclength, which will be called the \emph{density} of harmonic measure. 

The images of $\mathbb D$ under conformal contractions fixing $0$ are precisely those domains $\Omega$ for which the density of $\omega_\Omega(\cdot, 0)$ is at least $1/(2\pi)$ everywhere on the boundary. This follows immediately from the conformal invariance of harmonic measure and the fact that its density on the boundary of the unit disk is $1/(2\pi)$. Thus, Theorem~\ref{contra} gives a sufficient condition for $\Omega$ to have harmonic measure with such a lower density bound.

Theorem~\ref{contra} was prompted by a question of J. E. Tener~\cite{Te} which arose in the following context. When $f$ is a conformal map of $\mathbb D$ into itself with $f(0)=0$, the composition with $f$ is a contraction on the Hardy space $H^2(\mathbb D)$, see~\cite[Corollary 3.7]{CM}. By the conformal invariance of harmonic measure, this implies that the restriction operator $R\colon H^2(\mathbb D)\to L^2(\partial \Omega,  \omega_\Omega(\cdot, 0))$ is a contraction. A lower bound on the density of $\omega_\Omega(\cdot, 0)$ then allows one to estimate the norm of the restriction operator $R\colon H^2(\mathbb D)  \to  L^2(\partial \Omega)$  where $L^2$ is taken with respect to arclength. 

Concerning the structure of condition~\eqref{sufc}  it should be noted that the term
\[(R_O-R_C)\frac{\log R_I-\log R_C}{R_I-R_C}\]
is scale-invariant, while the second term, $\frac12 \log R_C$, tends to $-\infty $ as the domain is scaled down. Thus, for any convex domain $\Omega$ of class $C^{1,1}$  Theorem~\ref{contra} gives an explicit factor $\lambda>0$ such that the scaled-down domain $\lambda \Omega$ is the image of $\mathbb D$ under a conformal contraction. This can be compared to the classical Kellogg-Warschawski theorem~\cite[Theorem 3.5]{Po} which asserts that the conformal map of the disk onto a Dini-smooth Jordan domain $\Omega$ has a uniformly continuous derivative. The latter also implies that $\lambda \Omega$ is the image of $\mathbb D$ under a conformal contraction for sufficiently small $\lambda>0$. However, in contrast to Theorem~\ref{contra}, one does not have an explicit suitable value of  $\lambda$ in this case.

Examples illustrating and motivating the condition~\eqref{sufc} are given in ~\S\ref{secex}.

The higher-dimensional version of Theorem~\ref{contra} is stated in terms of the harmonic measure, since there is no longer a rich supply of conformal maps. The desired property of $\Omega$ in this case is having the density of $\omega_\Omega(\cdot, 0)$  at least $1/\sigma_{n-1}$, where $\sigma_{n-1}$ is the surface area of the unit sphere. The quantity $1/\sigma_{n-1}$ is the density of the harmonic measure of the unit ball with respect to its center. We will also use the notation
\[
a^+ = \max(a,0),\qquad \phi(a,b)=\frac{\log a-\log b}{a-b},\quad \phi(a,a)=1/a.
\]

\begin{theorem}\label{hithm} Let $\Omega\subset \mathbb R^n$, $n>2$, be a convex domain that contains $0$ and has $C^{1,1}$-smooth boundary. If the radii $R_O$, $R_I$, and $R_C$ satisfy 
\begin{equation}\label{hithm1}
R_C R_I^{n-2} e^{n (R_O-R_C-R_I/2)^+ \phi(R_I/2,R_C)} \le \frac{2^{n-2}-1}{2^{n-1}(n-2)}
\end{equation}
then the density of the harmonic measure of $\Omega$ with respect to $0$ is bounded below by $\sigma_{n-1}^{-1}$. 
\end{theorem}

The exponential term in ~\eqref{hithm1} is scale invariant, while the factor $R_CR_I^{n-2}$  makes sure that the left hand side of~\eqref{hithm1}  tends to $0$ as the domain is scaled down.

\section{Examples and counterexamples}\label{secex}

The sufficient conditions of Theorems~\ref{contra} and ~\ref{hithm} are not necessary; however, they are reasonably precise. For example, in the special case $R_O=R_I=R_C=R$ the hypothesis of Theorem~\ref{contra} is that $R\le 1$, which is both necessary and sufficient in this case. The higher-dimensional estimate is less accurate:  the inequality ~\eqref{hithm1} simplifies to 
$R\le \frac12((2^{n-2}-1)/(n-2))^{1/(n-1)}$, where the right hand side is less than $1$ but converges to $1$ as $n\to\infty$. 

To justify the presence of three radii $R_O$, $R_I$, $R_C$ in Theorems~\ref{contra} and ~\ref{hithm}, let us note that constraining just two of them would not be sufficient for the conclusion. Indeed, a convex polygon has zero density of harmonic measure at the vertices. Slightly rounding the corners, one obtains a domain that fails the conclusion of the theorem, which only the curvature radius detects.  To show the necessity of $R_I$, let $\Omega$ be the disk of radius $1$ centered at the point $1-\epsilon$; the density of $\omega_\Omega(\cdot,0)$ is small on most of the boundary. Finally, letting $\Omega$ be the convex hull of the union of two disks such as $D(0,1)\cup D(n,1)$ shows  that the presence of $R_O$  is also necessary.

\begin{figure}[h]\label{triangle}
\centering 
\includegraphics[width=0.3\textwidth]{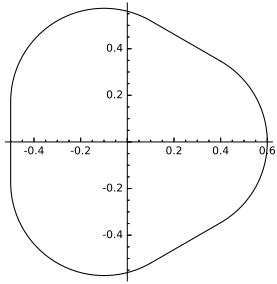}
\caption{A domain in Theorem \ref{contra}}
\end{figure}

Figure~1 presents a concrete example of a domain that satisfies~\eqref{sufc}, namely a rounded triangle with $R_O=0.6$, $R_I=0.5$, and $R_C=0.4$. 

A domain satisfying~\eqref{hithm1} could have $n=3$, $R_O=3/4$, and $R_I=R_C=1/2$. 

\section{Preliminaries: hyperbolic and quasihyperbolic metrics}\label{prelim} 

The hyperbolic metric on the unit disk $\mathbb D$ is 
\[
\rho_{\mathbb D}(z,w) = \inf\int_\gamma \frac{\abs{d\zeta}}{1-\abs{\zeta}^2} 
\]
where the infimum is taken over all rectifiable curves $\gamma$ connecting $z$ and $w$. In particular, 
\begin{equation}\label{hyper1}
\rho_{\mathbb D}(z,0) = \int_0^{\abs{z}} \frac{dt}{1-t^2}  = \frac12 \log\frac{1+\abs{z}}{1-\abs{z}}.
\end{equation}
On other simply-connected domains the hyperbolic metric can be defined by its conformal invariance property: 
$\rho_\Omega(f(z),f(w))= \rho_{\mathbb D}(z,w) $ if $f$ is a conformal map of $\mathbb D$ onto $\Omega$. 
In particular, for a disk $\Omega = D(a,R)$ we have 
\begin{equation}\label{hyper2}
\rho_{D(a,R)}(z,a) = \frac12 \log\frac{R+\abs{z-a}}{R-\abs{z-a}}.
\end{equation}

As a consequence of the Schwarz-Pick lemma, the  hyperbolic metric is  monotone with respect to domain: if $G$ and $\Omega$ are two simply-connected domains and $z,w\in G\subset  \Omega$, then 
\begin{equation}\label{hypermon}
\rho_{G}(z,w) \ge \rho_\Omega(z,w). 
\end{equation}

The quasihyperbolic metric $\rho^*_\Omega$ is defined by 
\[
\rho_{\Omega}^*(z,w) = \inf\int_\gamma \frac{\abs{d\zeta}}{\dist(\zeta,\partial \Omega)}. 
\]
It is not conformally invariant, but is comparable to $\rho_\Omega$ for every simply-connected domain: 
\begin{equation}\label{qh}
\frac14 \rho_\Omega^*(z,w) \le   \rho_\Omega(z,w) \le \rho_\Omega^*(z,w). 
\end{equation}
See~\cite[\S I.4]{GM} or ~\cite[\S 4.6]{Po}. 

\section{Planar domains: Proof of Theorem~\ref{contra}}\label{secplanar}

When the domain $\Omega$ is rescaled by the map $z\mapsto (1-\epsilon)z$, the left side of ~\eqref{sufc} decreases. Therefore, we may assume that strict inequality holds in~\eqref{sufc}. 

Let $f$ be a conformal map of the unit disk $\mathbb D$ onto $\Omega$, normalized by $f(0)=0$. By the Kellogg-Warschawski theorem~\cite[Theorem 3.5]{Po}, $f'$ has a continuous extension to $\overline{\mathbb D}$, and for $\zeta\in\partial \mathbb D$ we have 
\begin{equation}\label{der}
\lim_{z\to\zeta, \ z\in\mathbb D}\frac{f(z)-f(\zeta)}{z-\zeta} = f'(\zeta).
\end{equation}
By the maximum principle, it suffices to show $\abs{ f'}\le 1$ on $\partial\mathbb D$. By~\eqref{der} it suffices to show that 
\begin{equation}\label{der1}
\lim_{\abs{z}\nearrow 1 } \frac{\dist(f(z), \partial\Omega) }{1-\abs{z} }  \le 1.
\end{equation}

Fix $z\in\mathbb D$ and let $d=\dist(f(z), \partial\Omega)$. Since the small values of $d$ are of interest, we may assume $d<R_C$. Our plan is to estimate $\rho_\Omega(0,f(z))$ from above, which will yield 
\begin{equation}\label{con}
\abs{z} < 1-d 
\end{equation}
for  sufficiently small  $d$, thus proving~\eqref{der1}. 

Choose a point $w \in\partial \Omega$ such that $\abs{f(z)-w}=d$. By the definition of $R_C$, there is a disk $D = D(a , R_C)$ that has $w$ on its boundary and is contained in $\Omega$. Observe that $f(z)$ lies on the radius of this disk connecting $a$ to $w$, and therefore $\abs{f(z)-a} = R_C-d$. By~\eqref{hypermon} and~\eqref{hyper2}, 
\begin{equation}\label{pf1}
\rho_{\Omega}(f(z),a) \le \rho_{D(a,R_C)}(f(z),a) =  \frac12 \log\frac{2R_C-d}{d}
\le  \frac12 \log\frac{2R_C}{d}.
\end{equation}

To estimate $\rho_\Omega(a,0)$ we use the comparison with $\rho_\Omega^*$ stated in~\eqref{qh}. 
Since $\Omega$ contains $D(0,R_I)$ and $D(a,R_C)$, the convexity of $\Omega $ implies
\begin{equation}\label{pf2}
\dist(ta,\partial \Omega)\ge t R_C+(1-t)R_I, \quad 0\le t\le 1. 
\end{equation}
Integration along the line segment from $0$ to $a$ yields 
\begin{equation}\label{pf3}
\rho_\Omega^*(a,0) \le \abs{a} \int_0^1 \frac{dt}{t R_C+(1-t)R_I} 
= \abs{a} \phi(R_I, R_C). 
\end{equation}
Since $D(a,R_C)\subset \Omega\subset D(0,R_O)$, we have $\abs{a}\le R_O-R_C$. In conclusion, 
\begin{equation}\label{pf4}
\rho_\Omega(a,0) \le (R_O-R_C) \phi(R_I,R_C). 
\end{equation}

Suppose that ~\eqref{con} fails, that is, $\abs{z} \ge 1-d $. 
 From the conformal invariance of hyperbolic metric,  
\begin{equation}\label{pf5}
\rho_{\Omega}(f(z),0)  =    \frac12 \log\frac{1+\abs{z} }{1-\abs{z}} 
\ge  \frac12 \log\frac{ 2 - d }{d}. 
\end{equation}
Combining~\eqref{pf1}, \eqref{pf4}, and \eqref{pf5} we obtain 
\[
\frac12 \log\frac{2-  d }{d}   \le  \frac12 \log\frac{2R_C}{d} + (R_O-R_C) \phi(R_I,R_C), 
\]
hence 
\begin{equation}\label{con2}
\frac12 \log \left(1 -  \frac{d}{2}\right)  \le  \frac12 \log R_C + (R_O-R_C) \phi(R_I,R_C). 
\end{equation}
Since the right hand side of~\eqref{con2} is negative,  the inequality implies a lower bound on $d$. 
Therefore, \eqref{con} hold provided that $d$ is sufficiently small.  
This proves Theorem~\ref{contra}.

\section{Higher dimensions: proof of Theorem~\ref{hithm}}\label{sechigher}

In this section $\Omega$ is a convex domain in $\mathbb R^n$, $n>2$, and $0\in \Omega$.
The density of $\omega_\Omega(\cdot, 0)$ with respect to the surface measure of $\partial \Omega$ is related to Green's function $g_\Omega$ by 
\[
\omega_\Omega(E,0) = \int_E \frac{\partial g_\Omega }{\partial n}.
\]
Here the derivative is taken along the interior normal, and $g_\Omega$ is Green's function with pole at $0$, normalized by $g_\Omega(x,0)=\frac{1}{(n-2)\sigma_{n-1}}|x|^{2-n}+O(1)$ as $x\to 0$.

Thus, to prove that the density of harmonic measure is no less than $\sigma_{n-1}^{-1}$, it suffices to show that 
\begin{equation}\label{higoal}
g_\Omega(x)\ge  \frac{1+o(1)}{\sigma_{n-1}} \dist(x,\partial\Omega),\quad 
x\to\partial\Omega.
\end{equation}
To this end we use a lower bound for $g_\Omega$ in terms of the quasihyperbolic metric. An  estimate of this kind is given in Section 1.2 of~\cite{Ai}, namely $g_\Omega(x,0)\ge \exp(-A\rho_\Omega^*(x,0))$ with unspecified $A$. But we need a more explicit bound, since the presence of $A$ in the exponent does not allow one to conclude with~\eqref{higoal}. 

Fix $x\in \Omega $ and let $d=\dist(x, \partial\Omega)$, assuming $d<R_C$. Choose a point $w \in\partial \Omega$ such that $\abs{x-w}=d$. By the definition of $R_C$, there is a ball $B(a , R_C)$ that has $w$ on its boundary and is contained in $\Omega$. Since $\abs{x-a} = R_C-d$, Harnack's inequality \cite[Theorem 1.4.1]{AG} yields 
\begin{equation}\label{Hn}
g_\Omega(x) \ge \frac{(R_C-\abs{x-a})R_C^{n-2}}{(R_C+\abs{x-a})^{n-1}} g_\Omega(a) 
= \frac{d\, R_C^{n-2}}{(2R_C-d)^{n-1}} g_\Omega(a).  
\end{equation}

Since $B(0,R_I)\subset \Omega$, it follows that the restriction of $g_\Omega$ to $B(0,R_I)$ is minorized by Green's function of this ball: specifically,
\begin{equation}\label{minor}
g_\Omega(x) \ge \frac{1}{(n-2)\sigma_{n-1}}(|x|^{2-n}-R_I^{2-n}),\quad |x|<R_I.
\end{equation}
In particular, at the point $a'=\frac{R_I}{2}\frac{a}{|a|}$ we have 
\begin{equation}\label{greenball}
g_\Omega(a') \ge \frac{2^{n-2}-1}{(n-2)\sigma_{n-1}}R_I^{2-n}.
\end{equation}

As a corollary of Harnack's inequality \cite[Corollary 1.4.2]{AG}, the gradient of a positive harmonic function on $B(a,r)$ satisfies $|\nabla u(a)|\le (n/r) u(a)$.  Therefore, 
\[
|\nabla \log g_\Omega (x)|\le n/\dist(x,\partial\Omega')
\]
where $\Omega'=\Omega\setminus\{0\}$. This implies 
\begin{equation}\label{difflog}
|\log g_\Omega (x)-\log g_\Omega (y)|\le n \rho^*_{\Omega'}(x,y).
\end{equation}

\textit{Case 1:} $|a|\ge R_I/2$. Since $|a'| = R_I/2\le |a| \le R_O-R_C$ and $a'$ is a scalar multiple of $a$, it follows that $|a-a'| \le  (R_O-R_C-R_I/2)$. Observe that the domain $\Omega'$ contains the balls $B(a',R_I/2)$ and $B(a,R_C)$, as well as their convex hull. Integration similar 
to~\eqref{pf3} yields 
\[
\rho_{\Omega'}^*(a,a') \le |a-a'| \phi(R_I/2,R_C) \le (R_O-R_C-R_I/2) \phi(R_I/2,R_C).
\]
Using~\eqref{difflog} we obtain 
\[
g_\Omega(a) \ge \frac{2^{n-2}-1}{(n-2)\sigma_{n-1}}R_I^{2-n} e^{-n (R_O-R_C-R_I/2) \phi(R_I/2,R_C)}.
\]

\textit{Case 2:} $|a|<R_I/2$. Instead of using $a'$, we have 
\[
g_\Omega(a) \ge \frac{2^{n-2}-1}{(n-2)\sigma_{n-1}}R_I^{2-n}
\]
as in~\eqref{greenball}. 

Thus, in either case 
\[
g_\Omega(a) \ge \frac{2^{n-2}-1}{(n-2)\sigma_{n-1}}R_I^{2-n} e^{-n (R_O-R_C-R_I/2)^+ \phi(R_I/2,R_C)}
\]
which by virtue of~\eqref{Hn} implies   
\begin{equation}\label{hie1}
\frac{\sigma_{n-1}g_\Omega(x)}{d} \ge \frac{ R_C^{n-2}}{(2R_C-d)^{n-1}} \frac{2^{n-2}-1}{n-2 }R_I^{2-n} e^{-n (R_O-R_C-R_I/2)^+ \phi(R_I/2,R_C)}.
\end{equation}
As $d\to 0$, the right hand side of~\eqref{hie1} converges to 
\[
\frac{2^{n-2}-1}{2^{n-1}(n-2)}\frac{1}{R_CR_I^{n-2}} e^{-n (R_O-R_C-R_I/2)^+ \phi(R_I/2,R_C)} \ge 1.
\]
This proves~\eqref{higoal} and concludes the proof of Theorem~\ref{hithm}.  

\subsection*{Acknowledgement} I would like to thank the referee for carefully reading the paper and suggesting several improvements.

\end{document}